\documentclass{amsart}
\usepackage{enumerate}
\usepackage[dvips]{graphics}
\usepackage[dvips]{graphicx}
\usepackage{epsfig,enumerate,mathrsfs}
\usepackage{color}

\numberwithin{equation}{section}
\newtheorem{theorem}{Theorem}[section]

\newtheorem{definition}[theorem]{Definition}
\newtheorem{lemma}[theorem]{Lemma}
\newtheorem{corollary}[theorem]{Corollary}

\DeclareMathOperator{\dist}{dist}

\newcommand{\ppin}{\mathbb{P}_{\mathrm{in}}}
\newcommand{\ppout}{\mathbb{P}_{\mathrm{out}}}
\newcommand{\rin}{\mathrm{R}_{\mathrm{in}}}
\newcommand{\rout}{\mathrm{R}_{\mathrm{out}}}
\newcommand{\vterm}{\mathrm{V}_{\mathrm{term}}}
\newcommand{\vout}{\mathrm{V}_{\mathrm{out}}}
\newcommand{\vin}{\mathrm{V}_{\mathrm{in}}}
\newcommand{\tin}{\tau_{\mathrm{in}}}
\newcommand{\tout}{\tau_{\mathrm{out}}}

\def\N{\mathbb{N}}

\def\KD3{\mathrm{KD}^3}

\newcommand{\prob}{\mathbb{P}}

\def\max{\mathrm{max}}

\begin{document}

\title{Inverse problems for random walks on trees: network tomography}

\author{Victor de la Pena}
\address{Department of Statistics \\Columbia University\\
\\ NY, NY 13902\\ U.S.A. }
\email{vp@stat.columbia.edu}

\author{Henryk Gzyl}
\address{Department of Statistics \\University Simon Bolivar\\
\\ Caracas, Venezuela }
\email{hgzyl@usb.ve}

\author{Patrick McDonald}
\address{Division of Natural Science \\New College of Florida
\\Sarasota \\FL 34243 \\U.S.A.}
\email{mcdonald@ncf.edu}

\thanks{2000 Mathematics Subject Classification.
Primary: 60J10, 90B10}

\date{{\today. file name}: trees2.tex}

\begin{abstract}
Let $G$ be a finite tree with root $r$ and associate to the internal
vertices of $G$ a collection of transition probabilities for a simple
nondegenerate Markov chain.  Embedd $G$ into a graph $G^\prime$
constructed by gluing finite linear chains of length at least 2 to the
terminal vertices of $G.$  Then $G^\prime$ admits distinguished
boundary layers and the transition probabilities associated to the
internal vertices of $G$ can be augmented to define a simple
nondegenerate Markov chain $X$ on the vertices of $G^\prime.$ 
We show that the transition probabilities of $X$ can be recovered from
the joint distribution of first hitting time and first hitting place
of $X$ started at the root $r$ for the distinguished boundary layers  
of $G^\prime.$ 
\end{abstract}

\maketitle

%\tableofcontents

%%%%%%%%%%%%%%%%%%%%%%%%%%%%%%%%%%%%%%%%%%%%%%%%%%%%%%%%%%%%%%%%%
\section{Introduction}
%%%%%%%%%%%%%%%%%%%%%%%%%%%%%%%%%%%%%%%%%%%%%%%%%%%%%%%%%%%%%%%%%

Computing networks consist of hardware devices (hosts, routers, end
terminals, etc), together with a collection of connections between
pairs of such devices, along which packets of information are passed.
This rudimentary structure is readily modeled by a graph $\Lambda=(V,E)$
where the vertex set $V$ represents hardware devices and the edge set
$E\subset V\times V$ represents direct connections between devices.
If, in addition to the underlying graph structure, parameters are
associated to vertices and edges, it is possible to produce a more
accurate model of a given computing environment.  When  the parameters
defining the model (including the underlying graph) are dynamic,
network performance can be expected to vary, and network control
and/or predictability become issues of serious consequence.  For real
world applications, the first step in addressing such issues involves
the accurate monitering of the parameters which define the network.  

For networks modeled as above, it is often the case that
direct monitering of system parameters is impossible and one must
rely on inference methods to produce reliable estimates for parameter
values (cf \cite{CCLWY} for a recent survey).  The search for good
estimates provides a rich source of challenging inverse problems (for
related work see \cite{BDF}, \cite{LY}, \cite{RCW}, \cite{TYBW} and
references therein).  In this paper we formulate and solve one such
problem.   

The problems in which we are interested involve a fixed network topology
in which one designated device can send packets to a collection of
endusers and monitor packet arrival at enduser positions, but cannot
directly observe any behavior for devices between packet origin
and packet collection.  We are interested in using packet arrival
times to determine network parameters associated to devices which are
not directly observable.  Thus, we are interested in a type of {\em
  network tomography} problem.  To concisely state our results, we
begin by formalizing the discussion.  

Let $\Lambda=(V,E)$ be a finite tree with root vertex $r,$ terminal
vertices $\vterm(\Lambda)\subset V,$ and internal vertices $V\setminus
\vterm(\Lambda)$ (see section 2 for background and notation).  At each
terminal vertex of $\Lambda$ glue on a finite linear tree with vertex set of
size at least 2 (the size of the linear chain is determined by the
geometry of $\Lambda$ and the terminal vertex in question: see section 2).
Call the resulting graph, denoted $\Lambda^\prime,$ an {\it
  augmentation of $\Lambda$} and
note that $\Lambda$ is naturally embedded in $\Lambda^\prime.$  Note that since
each terminal vertex of $\Lambda$ has in effect been replaced by a linear
segment of length at least 2, it is possible to identify ``two layers
of boundary in $\Lambda^\prime.$''  More precisely, denote the outer boundary
layer of $\Lambda^\prime$ by $\vout(\Lambda^\prime) = \vterm(\Lambda^\prime)$
and denote the inner boundary layer of $\Lambda$ by
$\vin(\Lambda^\prime) = \{v^* \in V(\Lambda^\prime): v^*v \hbox{ is an
  edge}\}.$   

To each internal vertex of $\Lambda$ associate the transition probabilities
of a simple nondegenerate Markov chain.  Thus, for each internal
vertex $v^*$ there corresponds a collection $\{t_{v^*v}\}$ of positive
probabilities which sum to 1 as $v$ varies over vertices for which
$v^*v$ is an edge.  To every vertex of $\Lambda^\prime$ which is not an
internal vertex of $\Lambda$ or a terminal vertex of $\Lambda^\prime,$ associate
the transition probabilities of the simple symmetric random walk
(by construction, for every such vertex there are precisely two
adjacent vertices).  The transition probabilities associated to the
internal vertices of $\Lambda$ together with the transition probabilities
associated to the vertices of $\Lambda^\prime$ which are not internal
vertices of $\Lambda$ suffice to define a simple Markov chain on $\Lambda^\prime$
(killed upon reaching the outer boundary of $\Lambda^\prime$).  We denote
this Markov chain by $X$ and refer to it as the augmented chain
associated to $\Lambda^\prime.$  We ask:  
\begin{quote}What can we learn about the transition probabilities
  associated to the internal vertices of the network $\Lambda$ by monitoring
  the arrival times of the   Markov chain $X,$ started at the root of $\Lambda,$
  at detectors placed in the boundary layers?  
\end{quote}

Our main result is that there is an augmentation for which we can
recover {\em everything}.  Postponing technical definitions to section
2, we have  

\begin{theorem}\label{introtheorem} Let $\Lambda=(V,E)$ be a finite tree
  with root $r$ and unkown transition   probabilities associated to
  its internal vertices.  Let $\Lambda^\prime$ be the $2$-spherical augmentation
  of $\Lambda$ and suppose that $X$ is the augmented chain associated to
  $\Lambda^\prime.$  Let $\ppin^r$ be the joint distribution of the first hitting
  time and the first hitting place of $X$ started at $r$ for vertices
  in the inner boundary layer $\vin(\Lambda^\prime)$ and let $\ppout^r$ be the
  joint distribution of the first hitting time and the first hitting
  place of $X$ started at $r$ for vertices in the outer boundary layer 
  $\vout(\Lambda^\prime).$  Then $\ppin^r$ and $\ppout^r$ completely determine the
  transition probabilities of $X.$
\end{theorem}

Our proof of Theorem \ref{introtheorem} establishes somewhat more: we
show that the unkown transition probabilities are rational functions
of special values of the joint distributions $\ppin^r$ and $\ppout^r$
with explicit bounds on the time variable given in terms of the
geometry of the graph of $\Lambda$ (see Corollary \ref{algorithm}).  As a result,
the natural empirical statistics are consistent estimators for the
unkown transition probabilities.  Statistical extensions and
refinements of these results will appear elsewhere.

The proof of Theorem \ref{introtheorem} (see also Theorem
\ref{maintheorem}) relies on a close examination of the structure of
the path space associated to the process $X$ and recursion.  The
argument generalizes that given for chains in \cite{DGM1}; it depends
on the tree structure of $\Lambda.$  

While we have chosen to present our results in the context of
computing, it is clear that Theorem \ref{introtheorem} should have
applications in a variety of applied environments (our interests were
originally directed towards inverse problems for one-dimensional
diffusions).  Indeed, there is a variety of related literature (for
applications in medical imaging, see the survey \cite{A}; for an
application involving neuroscience, see \cite{BC}).  A number of such
applications have been discussed by Gr\"unbaum and his collaborators
(cf \cite{Gr1}, \cite{Gr2}, \cite{GM} \cite{P} and references therein).  Our 
results involve different techniques and focus on detailed ``time of
flight'' information, distinguishing it from the work cited above.

\section{Background and Notation}

Let $\Lambda=(V,E)$ be a finite rooted tree.  Thus, $\Lambda$ is a
connected graph without cycles, with finite vertex set $V$ and edge set $E\subset
V\times V,$ and a distinguished element $r \in V.$  We will say that a
vertex $v$ is {\em   terminal} if there is exactly one vertex $u$ such
that $uv\in E.$  We will write $\vterm(\Lambda)$ for the set of terminal
vertices.  We refer to vertices which are not terminal as {\em
  internal vertices.}   

By a path in $\Lambda$ we will mean an
ordered tuple of vertices, $(v_1,v_2, \dots, v_n)$ where $v_iv_{i+1}
\in E$ for all $i<n.$  Given a path $\gamma = (v_1,v_2, \dots, v_n)$
we will say that $\gamma$ connects $v_1$ to $v_n.$  We write
${\mathcal P}_{uv}$ for the collection of paths connecting $u$ and
  $v.$

Associated to every path $\gamma$ is a length: the length of $\gamma,$
denoted by $l(\gamma),$ is the number of edges defined by 
$\gamma$ (i.e. if $\gamma = (v_1,v_2, \dots, v_n),$ then $l(\gamma)
=n-1$).  There is a natural notion of distance between vertices of $\Lambda:$
\[
\dist(u,v) = \min_{\gamma \in {\mathcal P}_{uv}} \{l(\gamma)\}.
\]
The length function and the root give rise to a norm:
\[
|v| = \dist(v,r).
\]
The norm gives rise to a partition of the vertices of $\Lambda$ by
shells: 
\begin{equation}
V_k(\Lambda) = \{v\in V(\Lambda): |v|=k\}.\label{shells2.1}
\end{equation}
Because $\Lambda$ is a tree, if $v\in V_{k+1},$ there is a
{\em unique} vertex $v^*\in V_k$ such that $v^*v$ is an edge.

A simple but important example is given by a discrete interval: given
integers $-k\leq 0 <l,$ the $(-k,l)$-segment is the rooted tree 
with $k+l+1$ vertices obtained by taking as vertices the integers
$\{-k, -k+1,\dots,l\},$ and edges given by $i(i+1), \ -k\leq i<l,$ and
root vertex $0$ (when $k=0$ we will refer to the corresponding
segment as the $l$-segment).  For this example, the shells of
$\Lambda$ contain either one or two vertices.  

If $\Lambda$ is a finite tree with root $r,$ we can associate to
$\Lambda$ an outer radius and an inner radius:
\begin{eqnarray}
\rout(\Lambda) & = & \max_{v\in V(\Lambda)}\{|v|\} \\
\rin(\Lambda) & = & \min_{v\in \vterm(\Lambda)}\{|v|\}.
\end{eqnarray}
We will say that $\Lambda$ is spherical if $\rout(\Lambda) =
\rin(\Lambda).$  Note that for spherical rooted trees we have 
\[
V_{\rout(\Lambda)} = \vterm(\Lambda).
\]

\begin{definition}\label{augmentation}
Suppose that $\Lambda$ is a rooted tree, $v\in \vterm(\Lambda),$ and
  $l\in \N.$  The {\em $l$-augmentation of $\Lambda$ at $v$} is the
  finite rooted tree obtained by gluing a copy of the $l$ segment to
  $\Lambda$ by identifying the root of the $l$-segment with the
  terminal vertex $v$ and taking the root of the resulting tree to be
  the root of $\Lambda.$  We say that the rooted tree $\Lambda^\prime$
  is an augmentation of $\Lambda$ if $\Lambda^\prime$ is obtained from
  $\Lambda$ by a series of augmentations at boundary vertices.  We say
  that $\Lambda^\prime$ is an $m$-complete augmentation of $\Lambda$
  if $\Lambda^\prime$ is an augmentation of $\Lambda$ which can be
  obtained by first performing an $m$-augmentation of $\Lambda$ at
  each terminal vertex of $\Lambda.$  
\end{definition}

\begin{figure}
\begin{center}
\mbox{\epsfig{file=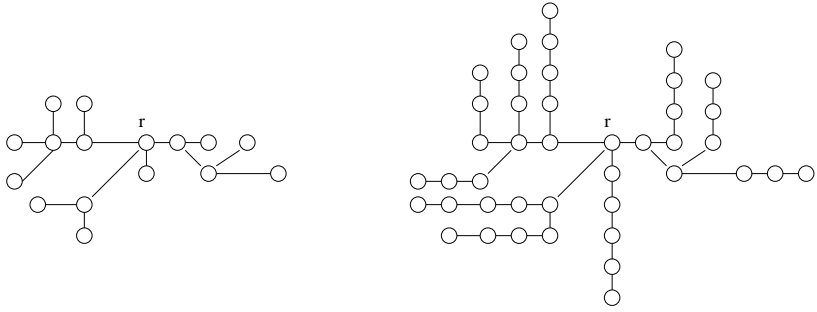}}
\end{center}
\caption{A rooted tree and a $2$-complete augmentation.}
\label{rgraph1}
\end{figure}

Thus, if we $l$-augment the $l$-interval at $0\in \vterm([0,l])$ we
obtain the $(-l,l)$-interval.  Similarly, given $n$ copies of the
$l$-interval, we can glue them together at $0$ and take as root the
gluing point to obtain a spherical rooted tree of radius $l$ which we
will refer to as the $(l,n)$-star.    

Given a rooted tree $\Lambda$ we can always construct an augmentation
of $\Lambda$ which is spherical:

\begin{definition}Suppose that $\Lambda$ is a rooted tree and that $l$
  is any natural number.  For each $v\in \vterm(\Lambda),$ perform an 
  $(\rout(\Lambda) -|v| + l)$-augmentation of $\Lambda$ at $v.$   The
  resulting augmentation is called the $l$-spherical augmentation
of $\Lambda.$  
\end{definition}

Thus, if $\Lambda $ is a finite rooted tree and $\Lambda^\prime$ is
  the $l$-spherical augmentation of $\Lambda,$ then $\Lambda^\prime$
  is spherical and  
\begin{equation*}
\rin(\Lambda^\prime) = \rout(\Lambda^\prime) = \rout(\Lambda)+l.
\end{equation*}

Given a finite rooted tree $\Lambda$ and an $l$-spherical augmentation
$\Lambda^\prime,$  suppose that $X$ is a nondegenerate simple Markov
chain on $\Lambda^\prime.$  We are interested in determining the
structure of $X$ on $\Lambda$ via first passage probabilities at 
the boundary of $\Lambda^\prime.$  More precisely,

\begin{definition}\label{determined}Let $\Lambda$ be a finite tree
  with root $r$ and
  $\Lambda^\prime$ an $l$-spherical augmentation of $\Lambda.$  Let
 $X$ be a nondegenerate simple Markov chain on $\Lambda^\prime.$  We
 say that {\em $X$ is determined by $l$-spherical first hitting times}
 if the transition  probabilities for vertices in $\Lambda$ are
 completely determined by  the triple $(W,\ppout^r,\ppin^r)$ where 
\begin{itemize}
\item  $W$ is the collection of transition   probabilities for
  vertices $V(\Lambda^\prime)\setminus V(\Lambda),$
\item $\ppout^r$ is the joint distribution of hitting time and hitting
  place of $X$ started at the root for vertices in the shell
  $V_{\rout(\Lambda^\prime)},$ and 
\item $\ppin^r$ the joint distribution of
  hitting time and hitting place  of $X$ started at the root for
  vertices in the shell $V_{\rout(\Lambda^\prime)-1}.$  
\end{itemize}
\end{definition}

We can now concisely state our main theorem:

\begin{theorem}\label{maintheorem} Let $\Lambda$ be a finite rooted
  tree and   suppose that $X$ is a nondegenerate simple Markov chain
  on the $2$-spherical augmentation of $\Lambda.$  Then $X$ is
  determined by $2$-spherical first hitting times. 
\end{theorem}

Before proceeding to the proof of Theorem \ref{maintheorem}, we note
that the main result of \cite{DGM1} establishes an important special
case: 

\begin{theorem}\label{dgm1result}\cite{DGM1} Let $k$ and $l$ be
  natural numbers and   suppose that $X$ is a simple nondegenerate
  Markov chain on the   $2$-spherical augmentation of the
  $(-k,l)$-interval.  Then $X$ is   determined by $2$-spherical first
  hitting times. 
\end{theorem}

Thus, as noted in \cite{DGM1}, there are simple
counterexamples to the most straightforward generalizations of Theorem
\ref{maintheorem} (in particular, it is, in general, impossible to
determine transition probabilities with a single pair of symmetrically
placed detectors, i.e. with a single boundary layer).  

The proof of Theorem \ref{maintheorem} involves an extension of the
ideas developed for the proof of Theorem \ref{dgm1result}.  As an
important illustrative  example of how the argument proceeds, we fix a
natural number $m$ and study the special case of the $(1,m)$-star. 

Enumerate the vertices of $(1,m)$-star $\Lambda$ as $\{v_i\}_{i=0}^m$
  with $v_0   =0.$  Using the norm to partition the vertices by
  shells (cf \eqref{shells2.1}), we write 
\begin{eqnarray*}
V_2 & = & \{v_{m+j}\}_{j=1}^m \hbox{ where there is an edge between
  $v_j$ and $v_{m+j}$}\\
V_3 & = &\{v_{2m+j}\}_{j=1}^m \hbox{ where there is an edge between $v_{m+j}$ and
  $v_{2m+j}.$}
\end{eqnarray*}

%\begin{figure}
%\begin{center}
%\mbox{\epsfig{file=pscat1.eps}}
%\end{center}
%\caption{Enumerating $V(A_2(T_{1,m})).$} 
%\label{onedimpicture}
%\end{figure}   

Let $\tin$ be the first hitting time of $V_2$
\begin{eqnarray*}
\tin & = & \inf\{n\geq 0 : X_n \in V_2 \}
\end{eqnarray*}
and similarly, let  $\tout$ be the first hitting time of $V_3.$  We will
write the joint distribution of first hitting time and first hitting
place as in Definition \ref{determined}:
\begin{eqnarray*}
\ppin^0(k,j) &= & \prob^0(\tin = k, X_{\tin} = v_{m+j}) \\ 
\ppout^0(k,j) &= & \prob^0(\tout = k, X_{\tout} = v_{2m+j}).
\end{eqnarray*}
We write $t_{j,l}$ for the probability of transitioning from vertex
$v_j$ to vertex $v_l$ in one time step.  Then 
\begin{eqnarray}
\ppin^0(2,j)& = & t_{0,j}t_{j,m+j} \label{2probs} \\
\ppout^0(3,j)& = & t_{0,j}t_{j,m+j}t_{m+j,2m+j}.
\end{eqnarray}
Similarly, 
\begin{eqnarray}
\ppin^0(4,j)& = & \ppin^0(2,j)\left[\sum_{i=1}^m
  t_{0,i}t_{i,0}\right] \\
\ppout^0(5,j)& = & \ppout^0(3,j)\left[\sum_{i=1}^m
  t_{0,i}t_{i,0} + t_{m+j,j}t_{j,m+j}\right].  
\end{eqnarray}
We conclude 
\begin{equation}
t_{j,m+j} = \frac{1}{t_{m+j,j}}\left[
\frac{\ppout^0(5,j)}{\ppout^0(3,j)} -
\frac{\ppin^0(4,j)}{\ppin^0(2,j)} \right]
\end{equation}
from which it follows that the transition probabilities $t_{j,m+j}$
are determined for all $j.$  From (\ref{2probs}), it follows that the
transition probabilities $t_{0,j}$ are determined for all $j.$  Since
the walk is by assumption simple, the transition probabilities
$t_{j,0}$ are determined for all $j.$  We conclude:

\begin{lemma}\label{reduction1}
For every $m,$ every nondegenerate simple Markov chain on the $(1,m)$
star is determined by $2$-spherical first hitting times.
\end{lemma}

\section{Proof of Theorem \ref{maintheorem}}

To establish the general case, we begin with an observation:  Every
rooted tree naturally embedds in its spherical augmentations.  Thus, 

\begin{lemma} Theorem \ref{maintheorem} is true for general rooted
  trees   if and only if it is true for all spherical rooted
  trees. 
\end{lemma}%\stopthm%

Let $\Lambda$ be a spherical rooted tree.  To prove Theorem
\ref{maintheorem} we give a careful analysis of the structure of paths
beginning at the root and having certain prescribed hitting properties
for the outer shells of $\Lambda^\prime$ where $\Lambda^\prime$ is a
(general) spherical augmentation of $\Lambda.$  We first demonstrate
that we can recursively determine transition probabilities  associated
to terminal vertices of $\Lambda.$  Using this result and the methods
employed to obtain it, we establish a recursion algorithm for
determining all unknown transition probabilities.  

To this end, let $R>0$ be the radius of the spherical rooted tree
$\Lambda.$ Let $m>0$ and let $\Lambda^\prime$ be the $(m+1)$-spherical
augmentation of $\Lambda.$  We refer to elements of $V_{R+m+1}$
(respectively, $V_{R+m}$ as outer boundary vertices (respectively,
inner boundary vertices).  Fix a vertex $v\in V_{R+m+1}$ and define
(for the remainder of the paper)    
\begin{eqnarray}
T & = & (R+m+1) + 2m. \label{defT}
\end{eqnarray}
We consider paths beginning at the root and having $T$ as the first
hitting time for the vertex $v\in V_{R+m+1}.$  More
precisely, we define 
\begin{eqnarray}
\Gamma_v & = & \{\gamma: \gamma(0) = r, \ \gamma(T) = v, \nonumber
\\
   &   & \hspace{1in} \ \gamma(j) \neq v \hbox{ for all }
j<T\}. \label{paths} 
\end{eqnarray}

Elements of $\Gamma_v$ do not visit most boundary vertices of
$\Lambda^\prime:$

\begin{lemma}\label{insmall} Let $v\in V_{R+m+1}$ and suppose $v^* \in
  V_{R+m}$ is the unique vertex such that $v^*v$ is an edge.  For
  $\Gamma_v$ as in (\ref{paths}), if $\gamma \in \Gamma_v,$ then
  $\gamma(j)\notin V_{R+m}\setminus \{v^*\} $ for all $j\leq T.$  
\end{lemma}

{\sc Proof} Let $\gamma$ be a curve in $\Gamma_v$ and suppose
$|\gamma(j)| = R+m,$ $\gamma(j) \neq v^*.$  Then $j > R+m-1.$  Denote 
by $\tau$ the first time that $\gamma$ visits $v.$  Then 
\begin{eqnarray*}
\tau & \geq & j + (m+1) + (m+2) \\
 & > & R+m-1 +3 + 2m.
\end{eqnarray*}
We conclude that $\tau >T,$ which completes the proof of the
lemma. %\stopthm%    

We partition $\Gamma_v$ by first hitting times of $v^*:$ 

\begin{lemma}\label{decomp} Let $v\in V_{R+m+1}$ be an element of the
  outer boundary layer of $\Lambda^\prime$ and let $\Gamma_v$ be as in
  (\ref{paths}).  Suppose that $v^*$ is the unique vertex for which
  $v^*v$ is an edge and define 
\begin{eqnarray}
\Gamma_{v,l} & = & \{\gamma \in \Gamma_v: \gamma(T-(2l-1)) = v^*, \nonumber
 \\
   &    & \hspace{1in} \ |\gamma(j)| <R+m \hbox{ for all } j <
 T-(2l-1)\}. \label{gk} 
\end{eqnarray}
Then
\begin{enumerate}
\item  $\Gamma_{v,i} \cap \Gamma_{v,j} = \emptyset$ if $i\neq j$ 
\item  $\cup_{l=1}^{m+1} \Gamma_{v,l} = \Gamma_v.$
\end{enumerate}
\end{lemma}

{\sc Proof}  If $i \neq j,$ elements of $\Gamma_{v,i}$ and $\Gamma_{v,j}$ have
different first hitting times of $v^*$ and thus $\Gamma_{v,i}\cap \Gamma_{v,j}
= \emptyset.$  Since every element of $\Gamma_v$ begins at the root and hits
$v^*$ by time $T,$ $\cup_{l=1}^{m+1} \Gamma_{v,l} = \Gamma_v.$    %\stopthm%

Initial segments of paths in $\Gamma_{v,l}$ define paths with nice first
hitting properties.  We make this precise:

\begin{lemma}\label{truncate} Let $1\leq l \leq m+1.$  For $\gamma\in
\Gamma_{v,l},$ the path ${\mathcal C}_l(\gamma)$ obtained by truncating $\gamma$ at
time $T-(2l-1)$  satisfies 
\begin{enumerate}
\item  ${\mathcal C}_l(\gamma)(0) = r$ 
\item  ${\mathcal C}_l(\gamma)(T-(2l-1)) = v^*$
\item  $| {\mathcal C}_l(\gamma)(j)| < R+m$ for all $j < T-(2l-1).$
\end{enumerate}
\end{lemma}

{\sc Proof} The first statement is obvious.  By definition, for
$\gamma\in \Gamma_{v,l},$ the first hitting time of $v^*$ is
$T-(2l-1).$  From this and Lemma \ref{insmall} we conclude that (2) and
(3) hold.    %\stopthm%

Truncation provides for a decomposition of paths in $\Gamma_{v,l}:$  Each
such path consists of an initial segment which has nice first hitting
properties, followed by an end segment which never visits the tree
$\Lambda.$  We make this precise:  Denoting starting positions by a
superscript, we write 
\begin{equation}\label{jhprobabilities}
\prob^r(t,v) = \hbox{Probability}(\{\hbox{starting at $r,$ $X$ first hits $v$
at time $t$}\}). 
\end{equation}
We have:
\begin{lemma}\label{factoring} Let $\Gamma_{v,l}$ be as defined in
  lemma \ref{decomp}.  For $1\leq l < m+1,$
\begin{eqnarray}
\prob^r(\Gamma_{v,l}) & = & \ppin^r(T-(2l-1),v^*)\chi_{v,l} \label{gk2}
\end{eqnarray}
where $\chi_{v,l}$ is an expression which involves only the transition
probabilities for vertices $z$ with $ R+1 \leq |z| < R+m+1.$
\end{lemma}

{\sc Proof}  By Lemma \ref{truncate}, each $\gamma \in \Gamma_{v,l}$ can be
decomposed as a path ${\mathcal C}_l(\gamma)$ starting at the root $r$ with first hitting
time of $\{v^*\}$ occurring at time $T-(2l-1),$ followed by a path of
length $(2l-1)$ which begins at $v^*$ and 
ends when it makes its first visit to $v.$  We will write 
\begin{eqnarray}
\tilde{\Gamma}_{v,l} & = & \{\gamma: \gamma(0) = v^*, \gamma(2l-1) =
v, \nonumber \\ 
  &    & \hspace{1in} \ |\gamma(j)| <R+m+1 \hbox{ for all } j <
2l-1\}.\label{ek} 
\end{eqnarray}
By choice of $l,$ if $\gamma\in \tilde{\Gamma}_{v,l},$ $|\gamma(j)| >R$ for
all $0 \leq j \leq 2l-1.$  Thus, if $\gamma \in \tilde{\Gamma}_{v,l}$ we
can compute $\prob^{v^*}(\{\gamma\}) $ in terms of the transition
probabilities associated to vertices $z$ with $R< |z| \leq R+m.$
Summing over all elements $\tilde{\Gamma}_{v,l}$ gives an expression 
\begin{eqnarray*}
\chi_{v,l} & = & \prob^{v^*}(\tilde{\Gamma}_{v,l})
\end{eqnarray*}
which involves only the transition probabilities for vertices $z$ with
$ R< |z| \leq R+m.$  Finally, using Lemma \ref{truncate}, we compute
$\prob^r(\Gamma_{v,l}):$ 
\begin{eqnarray*}
\prob^r(\Gamma_{v,l}) & = & \ppin^r(T-(2l-1), v^*) \chi_{v,l} 
\end{eqnarray*}
as required.   

The next result establishes that the transition probabilities at the
terminal vertices of $\Lambda$ are $m$-spherically determined by first
hitting times.  It is also the first step in an inductive proof of
Theorem \ref{maintheorem}.

\begin{lemma}\label{key}  Let $m>1.$  Let $\Lambda$ be a spherical
  tree with root $r$ and radius $R,$ $\Lambda^\prime$ the $(m+1)$-spherical
  augmentation of $\Lambda.$   Let $u\in V_{R},$ and suppose $w$ is
  the unique vertex of   $\Lambda^\prime$ for which $uw$ is an edge of
  $\Lambda^\prime.$   Then the transition probability $t_{uw}$ is
  determined by the triple $(W,\ppout^r,\ppin^r)$ where $W$ is the set of
  transition probabilities for vertices $V(\Lambda^\prime) \setminus
  V(\Lambda),$ $\ppout^r$ is the joint distribution of first hitting
  time and first hitting place for $V_{R+m+1},$ and $\ppin^r$ is the joint
  distribution of first hitting time and first hitting place for
  $V_{R+m}.$  
\end{lemma}

{\sc Proof} Let $v$ be a terminal vertex of $\Lambda^\prime$ a
distance $m+1$ from $u.$  Let $T$ be as in (\ref{defT}) and let $v^*$
be the unique vertex such that $v^*v$ is an edge of $\Lambda^\prime.$
From Lemma \ref{decomp} and Lemma \ref{factoring} we have  
\begin{eqnarray}
\prob^r(T, v) & = & \prob^r(\Gamma)  \nonumber \\
 & = & \sum_{l=1}^m \ppin^r(T-(2l-1), v^*) \chi_{v,l} +
 \prob^r(\Gamma_{v,m+1}).\label{decomp2}
\end{eqnarray}
We let $\gamma_*$ be the element of $\Gamma_{v,m+1}$ which ``changes
direction exactly twice.''  From (\ref{defT}) and (\ref{gk}), $\gamma_*$ is
the path which starts at $r,$ moves out radially $R+m$ units, moves in
radially $m$ units and moves out radially $m+1$ units.  Since we know
all transitions associated to $\gamma_*$ we can explicitly compute the
probability that $\gamma_*$ occurs:  
\begin{eqnarray}
\prob^r(\{\gamma_*\}) & = &  \chi_* t_{uw}  \label{gstar}
\end{eqnarray}
where 
\begin{eqnarray}
\chi_* & = & \ppout^r(R+m+1, v)\rho_{uv} \label{chi}
\end{eqnarray}
and $\rho_{uv}$ involves only transition probabilities along the
path of length $m+1$ from $u$ to $v$ (if we write the unique such path
as $u_0u_1\dots u_{m+1}$ with $u_0=u$ and $u_{m+1}=v,$ then $\rho_{uv}
= (\prod_{i=1}^{m-1} t_{u_{i+1}u_i})(\prod_{i=0}^{m-1}
t_{u_iu_{i+1}})$). 

Recall, an element of $\Gamma_{v,m+1}$ starts at position $r,$ first
hits position $v$ at time $R+m+1+2m,$ and first hits position
$v^*$ at time $R+m.$  Thus, if $\gamma \in \Gamma_{v,m+1}\setminus
\{\gamma_*\},$ then, as in Lemma \ref{factoring}, we may view $\gamma$
as a truncation followed by a path which never visits $u.$  Thus, as
in Lemma \ref{factoring} we can write
\begin{eqnarray}
\prob^r(\Gamma_{v,m+1} \setminus \{\gamma_*\}) & = & \ppin^r(R+m,v^*)
\chi_{v,m+1} \label{factormp2} 
\end{eqnarray}
where $\chi_{v,m+1}$ depends only on transition probabilities for
vertices $z$ such that $ R < |z| \leq R+m.$  Using (\ref{decomp2}),
(\ref{gstar}) and (\ref{factormp2}) we can solve for $t_{uw}:$
\begin{eqnarray}
t_{uw} & = & \frac{1}{\chi_*}\left[ \ppout^r(T, v) - \left[
\sum_{l=1}^{m+1} \ppin^r(T-(2l-1), v^*) \chi_{v,l} \right]\right].
\label{tlfinal} 
\end{eqnarray}
This complete the proof of the lemma.   %\stopthm%

The next result provides for the inductive step in the proof of
Theorem \ref{maintheorem}.  

\begin{lemma}\label{induct} Let $\Lambda $ be a rooted tree, $\Lambda^\prime$ the
  $2$-spherical augmentation of $\Lambda.$  Suppose that $X$ is a
  simple nondegenerate Markov chain on $\Lambda^\prime,$ and that the
  transition probabilities for the vertices $V(\Lambda^\prime)
  \setminus (\cup_{j=0}^k V_j)$ are known.  Let $u \in V_k$ and $w\in
  V_{k+1}$ be such that $uw$ is an edge.  Then the transition
  probability $t_{uw}$ is determined by the transition probabilities
  at the vertices $V(\Lambda^\prime) \setminus (\cup_{j=0}^k V_j)$ and
  the joint distributions of first hitting time and place, $\ppin^r$
  and $\ppout^r.$
\end{lemma}

{\sc Proof} Let $\vterm(w,\Lambda^\prime)$ be the terminal vertices of
$\Lambda^\prime$ which can be connected to $w$ by a path of length
$R+2 -(k+1).$  Let 
\begin{equation}\label{inner3.3}
  \vterm^*(w,\Lambda^\prime) = \{v^* \in V_{R+1}: \hbox{ there exists
  } v \in \vterm(w,\Lambda^\prime) \hbox{ with $v^*v$ an edge}\}.
\end{equation}
Set 
\begin{equation}\label{Tref2}
T_w = R+1 +(R+1-k) + (R+2-k) = 3R+4-2k
\end{equation}
and let 
\begin{eqnarray}
\Gamma_{uw} & = & \{\gamma: \gamma(0) = r, \ \gamma(T_w) \in
\vterm(w,\Lambda^\prime) , \nonumber \\
   &   & \hspace{1in} \ \gamma(j) \notin \vterm(w,\Lambda^\prime)
\hbox{ for all } j<T_w\}. \label{paths2} 
\end{eqnarray}
As in Lemma \ref{insmall}, if $\gamma \in \Gamma_{uw},$ then
$\gamma(j) \notin V_{R+1} \setminus \vterm^*(w,\Lambda^\prime)$ for
all $j<T_w.$  For $1\leq l \leq k+1,$ set 
\begin{eqnarray*}
\Gamma_{uw,l} & = & \{\gamma \in \Gamma_{uw}: \gamma(T_w-(2l-1)) \in
 \vterm^*(w,\Lambda^\prime) , \nonumber 
 \\
   &    & \hspace{1in} \ |\gamma(j)| <R+1 \hbox{ for all } j <
 T_w-(2l-1)\}. 
\end{eqnarray*}
Then, as in Lemma \ref{decomp}, the sets $\Gamma_{uw,l}$ partition
$\Gamma_{uw}.$ Moreover, as in Lemma \ref{truncate}, paths behave
nicely under truncation in that for all $l$ with $1\leq l \leq k+1,$
paths in $\Gamma_{uw,l}$ have truncations which start at the root $r$
and first hit $V_{R+1}$ in $\vterm^*(w,\Lambda^\prime)$ at time
$T_w-(2l-1).$  As in the case $k=R,$ if $1\leq l<k+1,$ a path $\gamma
\in \Gamma_{uw,l}$ can be decomposed as a path ${\mathcal C}_l(\gamma)$ which
starts at the root $r$ and first hits $V_{R+1}$ in
$\vterm^*(w,\Lambda^\prime)$ at time $T_w-(2l-1),$ followed by a path
that never visits the $k$ shell $V_k(\Lambda).$  As in Lemma \ref{factoring},
we conclude 
\begin{eqnarray}
\prob^r(\Gamma_{uw,l}) & = & \sum_{v^*\in \vterm^*(w,\Lambda^\prime)}
\ppin^r(T_w-(2l-1),v^*)\chi_{v^*,l} \label{gk22} 
\end{eqnarray}
where $\chi_{v,l}$ is an expression which involves only the transition
probabilities for vertices $z$ with $ k+1 \leq |z| < R+2.$  As in
Lemma \ref{key}, the paths $\Gamma_{uw,k+1}$ contain a distinguished
subset of elements: those with an initial segment which moves to a 
radial distance of $R+1$ in time $R+1,$ followed by a segment that
moves in a radial distance of $k$ units in $k$ time units, followed by
a segment which moves a radial distance of $k+1$ units (see
(\ref{Tref2})).  If we denote this subset by $\Gamma^*_{uw},$ then, as
in Lemma \ref{key}, we have  
\begin{equation}
\prob^r(\Gamma^*_{uw}) = t_{uw} \left[ \sum_{v\in \vterm(w,\Lambda^\prime)}
\ppout^r(R+2,v)\rho_{wv}\right]\label{almostdone1}
\end{equation}
where $\rho_{wv}$ involves only transition probabilities along the
path from $w$ to $v$ (and these transition probabilities are by
assumption known).  Finally, if $\gamma \in \Gamma_{uw,k+1} \setminus
\Gamma^*_{uw},$ then $\gamma(R+1) \in \vterm^*(w,\Lambda^\prime)$ and
$\gamma(R+1+2k) \in \vterm(w,\Lambda^\prime)$ which implies that
$\gamma$ does not visit the $k$ shell $V_k(\Lambda)$ once it has left it.  We
conclude 
\begin{equation}
\prob^r(\Gamma_{uw,k+1} \setminus \Gamma^*_{uw}) =   \sum_{v^*\in \vterm^*(w,\Lambda^\prime)}
\ppin^r(R+1,v^*)\chi_{wv^*}\label{almostdone2}
\end{equation}
where $\chi_{wv^*}$ involves only transition probabilities for
vertices $z$ satisfying $k<|z|\leq R+1.$  Using (\ref{gk22}),
\eqref{almostdone1}, and \eqref{almostdone2}, we can, as in Lemma
\ref{key}, solve for $t_{uw}.$  This concludes the proof of the lemma.

{\sc Proof of Theorem \ref{maintheorem}}  The proof is recursive; an
induction on distance to the inner boundary of the $2$-spherical
augmentation.  The formal argument is as follows:

By Lemma \ref{reduction1}, it suffice to consider the case of
spherical trees $\Lambda$ of radius $R$ where $R>0$ is arbitrary.  Let 
$\Lambda^\prime$ be the $2$-spherical augmentation of $\Lambda$ and
for $u\in V(\Lambda),$ let $d=d(u)$ be the distance of $u$ from the
inner boundary, $V_{R+1}(\Lambda^\prime),$ of $\Lambda^\prime.$  If $d=1,$ then
by Lemma \ref{key}, $t_{uw}$ is determined by $\ppin^r$ and
$\ppout^r.$  If the result holds when $d=k-1,$ by Lemma \ref{induct}
it is true when $d=k.$  This finishes the proof.

From the proofs of Lemma \ref{key} and Lemma \ref{induct}, we
note that, given a finite rooted tree, $\Lambda,$ embedded in its
$2$-spherical augmentation $\Lambda^\prime$ and a simple nondegenerate
Markov chain, we only require a finite number of values of the joint
distribution of exit time and place to determine the transition
probabilities for a simple nondegenerate Markov chain on the embedded  
tree.  More precisely, we have: 

\begin{corollary}\label{algorithm}  Let $\Lambda$ be a rooted tree
  with outer radius $\rout(\Lambda).$  Let $\Lambda^\prime$ be the
  $2$-spherical augmentation of $\Lambda$ and suppose that $X$ is a
  simple nondegenerate Markov chain of $\Lambda^\prime.$  Then there
  is an algorithmic procedure for explicitly determining the
  transition probabilities of $X$ on $\Lambda$ from the $2$-spherical
  hitting times.  The algorithm depends on data from the joint
  distribution of exit time and place for time $t\leq 3R+4.$  
\end{corollary} 

{\sc Proof} From the proof of Lemma \ref{induct}, to determine
transition probabilities for elements of the $k$ shell $V_k(\Lambda)$
we need to sample times up to $3R+4 -2k$ (cf \eqref{Tref2}).  The
corollary follows immediately.

\end{document}